\documentclass[12pt]{amsart}
\title[A note on Schr\"odinger equation]{A note on Schr\"odinger equation with linear potential and hitting times}
\author{Gerardo Hern\'andez-del-Valle}
\begin{document}
\maketitle
\begin{abstract}
In this note we derive a solution to the Schr\"odinger-type backward equation which satisfies a necessary boundary condition used in hitting-time problems [as described in Hern\'andez-del-Valle (2010a)]. We do so by using an idea introduced by  Bluman and Shtelen (1996) which is worked out  in Hern\'andez-del-Valle (2010b). This example is interesting since it is independent of the parameter $\lambda$,  namely:
\begin{eqnarray*}
\kappa(s,x)=\frac{x}{\sqrt{2\pi s}}\exp\left\{-\frac{(x+\int_0^sf'(u))^2}{2s}\right\}.
\end{eqnarray*}
and suggest a procedure to generating more vanishing solutions and $x=0$.
\end{abstract}
\section{The Example.}
In Hern\'andez-del-Valle (2010b) the author finds solutions to a Schr\"o\-dinger-type backward equation which satisfy some necessary boundary condition used in hitting-time problems [see Hernandez-del-Valle (2010a)]. Namely, the PDE of interest is
\begin{eqnarray}\label{schrodinger}
-\frac{\partial w}{\partial t}(t,x)+xf''(t)w(t,x)=\frac{1}{2}\frac{\partial^2w}{\partial x^2}(t,x)
\end{eqnarray}
which alternatively is related to the following expectation:
\begin{eqnarray*}
w(t,x)=\tilde{\mathbb{E}}\left[\exp\left\{-\int_t^sf''(u)\tilde{X}_udu\right\}\Big{|}\tilde{X}_t=x\right]
\end{eqnarray*}
where process $\tilde{X}$ is the so-called three-dimensional Bessel bridge and $f$ corresponds to a moving boundary. Furthermore it was shown in Hernandez-del-Valle (2010a) that $w$ should satisfy the following inequality
\begin{eqnarray}\label{inequality}
0\leq w(t,x)\leq h(s-t,x)\qquad \forall\, 0\leq t<s, x\geq 0,
\end{eqnarray}
where $h$ is the so-called derived heat source solution defined as
\begin{eqnarray*}
h(s,x):=\frac{x}{\sqrt{2\pi s^3}}\exp\left\{-\frac{x^2}{2s}\right\}.
\end{eqnarray*}
(Which is also the density of the first hitting time of one-dimensional standard Brownian motion to the fixed boundary $x$.)
Thus, for all time $t$ it follows that $w$ should satisfy the following boundary condition:
\begin{eqnarray}\label{boundary}
\lim\limits_{x\to 0} w(t,x)=0.
\end{eqnarray}
Yet, for the specific problem of finding the density of hitting a moving boundary we only need the previous inequality (\ref{inequality}) and boundary condition (\ref{boundary}) to hold at $t=0$:
\begin{eqnarray}\label{atzero}
w(0,x)=\tilde{\mathbb{E}}\left[\exp\left\{-\int_0^sf''(u)\tilde{X}_udu\right\}\Big{|}\tilde{X}_0=x\right].
\end{eqnarray}
This is accomplished by using an idea introduced in Bluman and Shtelen (1996)
which relates the following PDE
\begin{eqnarray*}\label{backward}
-\frac{\partial u}{\partial t}(t,x)+V_1(t,x)u(t,x)=\frac{1}{2}\frac{\partial^2u}{\partial x^2}(t,x)
\end{eqnarray*}
and its adjoint
\begin{eqnarray}\label{forward}
\frac{\partial \Phi}{\partial t}(t,x)+V_1(t,x)\Phi(t,x)=\frac{1}{2}\frac{\partial^2\Phi}{\partial x^2}(t,x)
\end{eqnarray}
with the following backward equation
\begin{eqnarray*}
-\frac{\partial w}{\partial t}(t,x)+V_2(t,x)w(t,x)=\frac{1}{2}\frac{\partial^2w}{\partial x^2}(t,x)
\end{eqnarray*}
where
\begin{eqnarray*}
V_2(t,x)=V_1(t,x)-\frac{\partial^2}{\partial x^2}\log \Phi.
\end{eqnarray*}
It is done so through the following expression:
\begin{eqnarray}\label{new1}
w(t,x)=\frac{1}{\Phi(t,x)}\left[\int_0^xu(t,\xi)\Phi(t,\xi)d\xi+B_2(t)\right]
\end{eqnarray}
with $B_2(t)$ satisfying the condition
\begin{eqnarray}\label{new2}
\frac{dB_2}{dt}=\frac{1}{2}\left(\frac{\partial\Phi}{\partial x}(t,0)u(t,0)-\Phi(t,0)\frac{\partial u}{\partial x}(t,0)\right).
\end{eqnarray}
Thus at $t=0$, see equation (\ref{atzero}),
\begin{eqnarray*}
w(0,x)=\frac{1}{\Phi(0,x)}\left[\int_0^xu(0,\xi)\Phi(0,\xi)d\xi\right]
\end{eqnarray*}
and hence 
$$\lim\limits_{x\to0}w(0,x)=0$$
as long as 
$$\lim\limits_{x\to0}\Phi(0,x)\not=0.$$
Of course, in general we will not be solving the same PDE unless
\begin{eqnarray}\label{condition}
\frac{\partial^2}{\partial x^2}\log\Phi=0,
\end{eqnarray}
in which case $V_2(t,x)=V_1(t,x)$, that is the form of the PDE is preserved. The reader may check [or consult at Hernandez-del-Valle (2010a) and (2010b)] that given $V_1(t,x)=xf''(t)$ a solution to (\ref{forward}) which also satisfies condition (\ref{condition}) is for instance:
\begin{eqnarray*}
&&\Phi(t,x)\\
&&\quad=\exp\left\{\frac{1}{2}\int_0^t(f'(u))^2du-xf'(t)-\frac{1}{2}\lambda^2t-i\lambda\left(x-\int_0^tf'(u)du\right)\right\}
\end{eqnarray*}
where $i=\sqrt{-1}$ and $\lambda$ is some scalar. 

In the remainder of this note we derive a solution to (\ref{schrodinger}), which is independent of $\lambda$, and also satisfies boundary condition (\ref{boundary}) at $t=0$. We do so by using  Bluman and Shtelen's representation, equations (\ref{new1}) and (\ref{new2}), and a solution to (\ref{backward}) given by:
\begin{eqnarray}\label{u}
u(t,x)&=&\exp\left\{\frac{1}{2}\int_t^s(f'(u))^2du+xf'(t)\right\}\\
&&\qquad\times\exp\left\{-\frac{1}{2}\lambda^2(s-t)+i\lambda\left(x+\int_t^sf'(u)du\right)\right\}.
\end{eqnarray}
[The reader may consult Hern\'andez-del-Valle (2007).] It follows that:
\begin{eqnarray*}
\Phi(t,x)u(t,x)=\exp\left\{\frac{1}{2}\int_0^s(f'(u))^2du-\frac{1}{2}\lambda^2s+i\lambda \int_0^sf'(u)du\right\}
\end{eqnarray*}
and
\begin{eqnarray*}
\frac{dB_2}{dt}(t)&=&\frac{1}{2}\Phi u\left(\frac{\Phi_x}{\Phi}-\frac{u_x}{u}\right)\\
&=&-(f'(t)+i\lambda)u\Phi.
\end{eqnarray*}
Hence, $B_2$ might be written in the two following ways
\begin{eqnarray}\label{b1}
B_2(t)=-\left(\int_0^tf'(u)du+i\lambda t\right)u\Phi
\end{eqnarray}
or
\begin{eqnarray}\label{b2}
B'_2(t)=\left(\int_t^sf'(u)du+i\lambda(s-t)\right)u\Phi.
\end{eqnarray}
From (\ref{new1}) $w$ satisfies
\begin{eqnarray*}
w(t,x)&=&\frac{1}{\Phi}\left[\int_0^xu\Phi dy+B_2(t)\right]\\
&=&x\frac{u\Phi}{\Phi}+\frac{B_2(t)}{\Phi},
\end{eqnarray*}
and hence for $B_2$ as in  (\ref{b1})
\begin{eqnarray}\label{w1}
w=\left(\left\{x-\int_0^tf'(u)du\right\}-i\lambda t\right)u,
\end{eqnarray}
and for $B'_2$ as in  (\ref{b2})
\begin{eqnarray}\label{w2}
w=\left(\left\{x+\int_t^sf'(u)du\right\}+i\lambda(s-t)\right)u.
\end{eqnarray}

Finally, recall the following Fourier representations:
\begin{eqnarray}
\label{heat}k(t,x)&:=&\frac{1}{\sqrt{2\pi t}}e^{-\frac{x^2}{2t}}\\
\nonumber&=&\frac{1}{2\pi}\int_{-\infty}^{+\infty}e^{-\frac{1}{2}\lambda^2t+i\lambda x}d\lambda\\
\label{heatderived} h(t,x)&:=&\frac{x}{\sqrt{2\pi t^3}}e^{-\frac{x^2}{2t}}\\
\nonumber&=&\frac{1}{2\pi}\int_{-\infty}^{+\infty}(-i\lambda)e^{-\frac{1}{2}\lambda^2t+i\lambda x}d\lambda
\end{eqnarray}
also known as the source and {\it derived} source heat equations respectively.

After contour integration with respect to $\lambda$, (\ref{w1}) becomes:
\begin{eqnarray*}
w(t,x)&=&\exp\left\{\frac{1}{2}\int_t^s(f'(u))^2du+xf'(t)\right\}\\
&&\times \Bigg{\{}\left(x-\int_0^tf'(u)du\right)k\left(s-t,x+\int_t^sf'(u)du\right)\\
&&\quad+t\cdot\frac{\left(x+\int_t^sf'(u)du\right)}{(s-t)}k\left(s-t,x+\int_t^sf'(u)du\right)\Bigg{\}}
\end{eqnarray*}
and for $B'_2$ or equation (\ref{w2})
\begin{eqnarray*}
w(t,x)&=&\exp\left\{\frac{1}{2}\int_t^s(f'(u))^2du+xf'(t)\right\}\\
&&\times\Bigg{\{}\left(x+\int_t^sf'(u)du\right)k\left(s-t,x+\int_t^sf'(u)du\right)\\
&&\quad- \left(x+\int_t^sf'(u)du\right)k\left(s-t,x+\int_t^sf'(u)du\right)\Bigg{\}}\\
&\equiv&0.
\end{eqnarray*}

The second solution is identically zero, and the first evaluated at $t=0$ is
\begin{eqnarray*}
w(0,x)=\exp\left\{\frac{1}{2}\int_0^s(f'(u))^2du+xf'(0)\right\}x\cdot k\left(s,x+\int_0^sf'(u)du\right)
\end{eqnarray*}
This alternatively implies that an approximation to the first hitting time density is given by:
\begin{eqnarray*}
\frac{x}{\sqrt{2\pi s}}\exp\left\{-\frac{(x+\int_0^sf'(u)du)^2}{2s}\right\}.
\end{eqnarray*}

\subsection{More vanishing solutions at $x=0$ and $t=0$.} Observe that if $u$ is as in (\ref{u}) which alternatively solves (\ref{backward}) then
$$u'(t,x)=\Gamma(\lambda,s)u(t,x)$$
is also a solution to (\ref{backward}) for an arbitrary function $\Gamma$. For instance, suppose that $\Gamma(\lambda,s)=(-i\lambda)$ then equation (\ref{w1}) becomes:
\begin{eqnarray*}
w'&=&\left(\left\{x-\int_0^tf'(u)du\right\}-i\lambda t\right)(-i\lambda)u\\
&=&\left(\left\{x-\int_0^tf'(u)du\right\}(-i\lambda)+(-i\lambda)^2 t\right)u.
\end{eqnarray*}
After contour integration and observing that
\begin{eqnarray*}
\frac{1}{2\pi}\int_{-\infty}^{+\infty}(-i\lambda)^2e^{-\frac{1}{2}\lambda^2 t+i\lambda x}d\lambda=\left(\frac{x^2}{t^2}-\frac{1}{t}\right)k(t,x)
\end{eqnarray*}
it follows that the new solution $w'$ is given by
\begin{eqnarray*}
&&w'(t,x)\\
&&=\exp\left\{\frac{1}{2}\int_t^s(f'(u))^2du+xf'(t)\right\}\\
&&\enskip\times \Bigg{\{}\left(x-\int_0^tf'(u)du\right)h\left(s-t,x+\int_t^sf'(u)du\right)\\
&&\quad+t\cdot\left[\frac{\left(x+\int_t^sf'(u)du\right)^2}{(s-t)^2}-\frac{1}{(s-t)}\right]k\left(s-t,x+\int_t^sf'(u)du\right)\Bigg{\}}
\end{eqnarray*}
or
\begin{eqnarray*}
w'(0,x)=\exp\left\{\frac{1}{2}\int_0^s(f'(u))^2du+xf'(0)\right\}x\cdot h\left(s,x+\int_0^sf'(u)du\right)
\end{eqnarray*}
and $h$ is as in (\ref{heatderived}).

\end{document}